\documentclass[12pt,twoside]{article}
\usepackage{graphicx}
\usepackage{amsfonts,amsbsy,amssymb,amsmath}
\allowdisplaybreaks
\usepackage{cite}
\usepackage{graphics}
\usepackage{xcolor}
\usepackage{subfigure}
\usepackage{enumerate}
\def\bpsp{\begin{pspicture}}
\def\epsp{\end{pspicture}}

\newtheorem{theorem}{Theorem}[section]
\newtheorem{remark}[theorem]{Remark}
\newtheorem{example}[theorem]{Example}
\newtheorem{lemma}[theorem]{Lemma}
\newtheorem{corollary}[theorem]{Corollary}
\newtheorem{definition}[theorem]{Definition}
\newtheorem{proposition}[theorem]{Proposition}

\newtheorem{note}{Note}
\newtheorem{case}{Case}

\newtheorem{conjecture}{Conjecture}
\newtheorem{question}{Question}

\newcommand{\bea}{\begin{eqnarray}}
\newcommand{\eea}{\end{eqnarray}}
\newcommand{\beq}{\begin{eqnarray*}}
\newcommand{\eeq}{\end{eqnarray*}}

\def\m4{\mbox{\rm ~(mod $4$)}}

\def \bd{\begin{definition}}
\def \ed{\end{definition}}
\def \bqu{\begin{question}}
\def \equ{\end{question}}
\def \bcc{\begin{conjecture}}
\def \ecc{\end{conjecture}}
\def \bt{\begin{theorem}}
\def \et{\end{theorem}}
\def \bl{\begin{lemma}}
\def \el{\end{lemma}}
\def \bc{\begin{corollary}}
\def \ec{\end{corollary}}
\def \be{\begin{equation}}
\def \ee{\end{equation}}
\def \ben{\begin{enumerate}}
\def \een{\end{enumerate}}
\def \ba{\begin{array}}
\def \ea{\end{array}}
\def \bp{\begin{proposition}}
\def \ep{\end{proposition}}
\def \bx{\begin{example}}
\def \ex{\end{example}}
\def \br{\begin{remark}}
\def \er{\end{remark}}
\def \bdsc{\begin{description}}
\def \edsc{\end{description}}

\def \bn{\begin{case}}
\def \en{\end{case}}
\def \bnt{\begin{note}}
\def \ent{\end{note}}
\def\1{1\!\!1}

\def\mm2{\mbox{\rm ~(mod $2$)}}
\def\m4{\mbox{\rm ~(mod $4$)}}

\def\qed{\nolinebreak\hfill\rule{.2cm}{.2cm}\par\addvspace{.5cm}}

\def\m{\mu}

\def\1{\textbf{1}}
\def\0{\textbf{0}}

\begin{document}
\title{On maximum spectral radius of $\{H(3,3),~H(4,3)\}$-free graphs}
\author{Amir Rehman$ ^{a}$, S. Pirzada$ ^{b} $ \\
$^{a,b}${\em Department of Mathematics, University of Kashmir, Srinagar, India}\\
$^{a}$\texttt{aamirnajar786@gmail.com}; $ ^{b} $\texttt{pirzadasd@kashmiruniversity.ac.in}\\
}

\date{}

\pagestyle{myheadings} \markboth{Rehman, Pirzada}{On maximum spectral radius of $\{H(3,3),~H(4,3)\}$-free graphs}
\maketitle
\vskip 5mm
\noindent{\footnotesize \bf Abstract.} Let $G$ be a simple connected graph of size $m$. Let $A$ be the adjacency matrix of $G$ and let $\rho(G)$ be the spectral radius of $G$. A graph is said to be $H$-free if it does not contain a subgraph isomorphic to $H$. Let $H(\ell,3)$ be the graph formed by taking a cycle of length $\ell$ and a triangle on a common vertex. Recently, Li,  Lu and Peng [Y. Li, L. Lu, Y. Peng, Spectral extremal graphs for the bowtie, Discrete Math. 346(12) (2023) 113680.] showed that the unique $m$-edge $H(3,3)$-free spectral extremal graph is the join of $K_2$ with an independent set of $\frac{m-1}{2}$ vertices if $m\ge 8$ and the condition $m\ge 8$ is tight. In particular, if $G$ does not contain $H(3,3)$ as induced subgraph, they proved that $\rho(G) \leq  \frac{1+\sqrt{4m-3}}{2} $ and equality holds when $G$  is isomorphic to $S_{\frac{m+3}{2},2}$. Note that Li et al. denoted $H(3,3)$ by $F_2$. In this paper, we find the maximum spectral radius and identify the graph with the largest spectral radius among all \{$H(3,3), H(4,3)$\}-free graphs of size odd $m$, where $m\geq 259$. Coincidentally, we show that $\rho(G) \leq  \frac{1+\sqrt{4m-3}}{2}$ when $G$ forbids both $H(3,3)$ and $H(4,3)$. In our case, the equality holds when $G$ is isomorphic to the same graph.
\vskip 3mm



\section{Introduction}
\begin{figure}[h]
\centering
\subfigure[]{\includegraphics[width=0.20\linewidth]{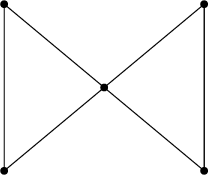}}~~~~~~~~~~
\subfigure[]{\includegraphics[width=0.24\linewidth]{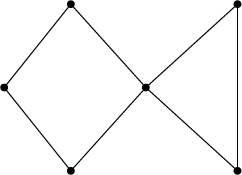}}
\caption{Graphs (a) H(3,3) and (b)  H(4,3)}
\label{fig01}
\end{figure}

 \indent Let  $G$ be a simple graph with order $n$ and size $m$ and let $V(G)$ be the vertex set of $G$. The adjacency matrix of $G$ is defined as $A(G)=(a_{ij})$, where
 \[
a_{ij}= \begin{cases}
      1 & \text{if there is an edge between vertices } i \text{ and } j, \\
      0 & \text{otherwise.}
\end{cases}
\]
\indent The largest eigenvalue of $A(G)$, denoted by $\rho$, is called the spectral radius of $G$. In case of a connected graph $G$, the Perron-Frobenius theorem asserts the existence of a unique positive eigenvector associated with $\rho(G)$, termed as the Perron vector of $G$. Additional definitions and notations can be found in \cite{CRS,SP}.\\
\indent For a subset $S\subseteq V(G),~G[S]$ represents the subgraph of $G$ induced by $S$. Further $e(S,T)$ denotes the number of edges with one end in $S$ and the other in $T$, where $S$ and $T$ are subsets of $V(G)$. Also, we use $e(S)$ to denote $e(S,S)$. We write $N^k(v)$ for the set of vertices at a distance of $k$ from vertex $v$, with $N^1(v)$ being denoted by  $N(v)$. We define $N[v]$ as $N(v)\cup\{v\}$. For $S \subseteq V(G)$, let $N_S(v)$ represent the set of neighbors of $v$ in $S$ and $d_S(v)$ be the cardinality of $N_S(v)$.\\
\indent For $1 \leq k \leq n$, the graph $S_{n,k}$ of order $n$ is obtained by joining each vertex of the complete graph $K_k$ to $n-k$ isolated vertices. Let $C_n$ represent the cycle on $n$ vertices. Let $H(\ell,3)$ be the graph formed by a cycle of length $\ell$ and a triangle on a common vertex. For example, the graphs $H(3,3)$ and $H(4,3)$ are shown in Figure \ref{fig01}. Define $G(m,t)$ to be the graph of size $m$ obtained by joining a vertex of maximum degree in $S_{\frac{m-t+3}{2},2}$ to $t$ isolated vertices (see Fig. \ref{fig 6}(a)). For $t=0$, the graph  $S_{\frac{m+3}{2},2}$ is called the book graph. Let $K^m_4$ be the graph of size $m$ obtained by joining a vertex from $K_4$ to $m-6$ isolated vertices (see Fig. \ref{fig 6}(b)).\\
\indent For a family of graphs $\mathcal{H}$, a graph $G$ is said to be $\mathcal{H}$-free if it does not contain a subgraph isomorphic to any graph in $\mathcal{H}$. In particular, if $\mathcal{H}=\{H\}$, we simply say that $G$ is $H$-free. A classical problem in the extremal graph theory is the Tur\'an problem which asks for the maximum size of an $H$-free graph of order $n$, where the maximum size is known as the Tur\'an number of $H$. Nikiforov \cite{NF3}, proposed a spectral analogue of the Tur\'an problem which asks for the maximum spectral radius of an $H$-free graph of size $m$ or order $n$.  In \cite{EN}, Nosal proved that $\rho(G) \le \sqrt{m}$ for every graph of size $m$, when $H$ is a triangle. Nikiforov \cite{NF2} showed that $\rho(G) \le \sqrt{m}$ for all $C_4$-free graphs of size $m$. Zhai, Lin, Shu \cite{ZHAI} showed that $\rho(G)\le \frac{1+\sqrt{4m-3}}{2}$ for any $C_5$-free graph of size $m\ge8$ or $C_6$-free graph of size $m\ge22$, with equality if and only if $G\cong S_{\frac{m+3}{2},2}$. In \cite{LP}, Li and Peng determined the maximum spectral radius of graphs with no intersecting odd cycles. Let $F_k$ be the graph obtained from $k$ triangles sharing a common vertex. Then $F_2$ is the graph $H(3,3)$. In 2023, Li, Lu and Peng \cite{LLP} characterized $F_2$-free graphs with given number of edges. They proved that the unique $m$-edge $H(3,3)$-free spectral extremal graph is the join of $K_2$ with an independent set of $\frac{m-1}{2}$ vertices (that is, $ S_{\frac{m+3}{2},2}$ ) if $m\ge 8$, and the condition $m\ge 8$ is tight. The problem has been investigated for various graphs $H$, as can be seen in \cite{CLZ1,CLZ,CDT,FYZ,NF,NF1,NF2,ZWF,ZW}.\\

The objective of this paper is to determine the maximum spectral radius of the graphs when  $\mathcal{H}= \{H(3,3), H(3,4)\}$. This is stated in the following theorem.
\begin{theorem}
          If $G$ is an \{$H(3,3),H(3,4)$\}-free graph with odd size $m\ge259$, and $G$ contains no isolated vertices, then $ \rho(G)\leq \frac{1+\sqrt{4m-3}}{2}$, unless $G\cong S_{\frac{m+3}{2},2}.$
\end{theorem}
\indent Coincidentally, the above bound is same as obtained by Li et al. \cite{LLP} for $H(3,3)$-free graphs.

 \indent In this paper, we use $G^\ast$ to represent the connected graph with maximum spectral radius among all graphs of size $m$ that are free of induced subgraphs $H(3,3)$ and $H(4,3)$. Let $\rho^\ast = \rho(G^\ast)$, and consider $X^\ast$ as the Perron vector of $G^\ast$ with coordinates $x_u$ corresponding to the vertex $u \in V(G^\ast)$. Let $x_{u^\ast} = \max\{x_u : u \in V(G^\ast)\}$, representing the coordinate associated with the vertex $u^\ast$ in $G^\ast$. We term $G^\ast$ as the extremal graph and $u^\ast$ as the extremal vertex of $G^\ast$.\\
Define $A(G^\ast) = A$, and let $N_{0}(u^\ast) = \{v : v \in N(u^\ast), d_{N(u^\ast)}(v) = 0\}$, and $N_{1}(u^\ast) = N(u^\ast)\setminus N_{0}(u^\ast)$. Additionally, define $N^2_{j}(u^\ast) =\{w\in N^2(u^\ast):d_{ N_{j}(u^\ast)}(w)\ge1\}$ for $j=0,1$ and let $W=V(G)\setminus N[u^\ast]$. By utilizing the eigenequations of $A$ on $u^\ast$, the following equation is obtained
\begin{equation}\label{eqn(21)}
\rho^\ast x_{u^\ast} = (AX^\ast){u^\ast} = \sum_{u \in N_0(u^\ast)} x_u + \sum_{v \in N(u^\ast)\setminus N_0(u^\ast)} x_v.
\end{equation}
As ${\rho^\ast}^2$ represents the spectral radius of $A^2$, the eigenequations of $A^2$ on $u^\ast$ lead to

\begin{equation}\label{eqn(22)}
{\rho^\ast}^2x_{u^\ast} = d(u^\ast)x_{u^\ast} + \sum_{u \in N(u^\ast)\setminus N_0(u^\ast)} d_{N(u^\ast)}(u)x_u + \sum_{w \in N^2(u^\ast)} d_{N(u^\ast)}(w)x_w.
\end{equation}

\indent The rest of the paper is organized as follows. In Section 2, we present lemmas which will be required to prove Theorem 1.1 in Section 3.
\section {\textbf{Lemmas}}

\indent In this section, we present several lemmas that will be required to prove the main theorem.
\begin{lemma}{\em \cite{NF2}} \label{lem21}
Let $A $ and $A^\prime$ be the adjacency matrices of two graphs $G$ and $G^\prime$ on the same vertex set. Suppose that $N_G(u) \subsetneq N_{G^\prime}(u)$ for some vertex $u$. If some positive eigenvector to $\rho(G)$ satisfies $X^\prime A^\prime X \ge X^\prime AX$, then $\rho(G^\prime) \textgreater \rho(G)$.
\end{lemma}
\begin{lemma}{\em \cite{CRS}}\label{lembap}
Let $G$ be a connected graph and let $H$ be a proper subgraph of $G$. Then $\rho(H)~\textless~\rho(G)$.
\end{lemma}
\begin{definition}{\em \cite{CRS}}
Given a graph $G$, the vertex partition $P:V(G)=V_1\cup V_2\cup \cdots \cup V_k$ is said to be an equitable partition if, for each $v\in V_i, ~|V_j\cap N(v)|=c_{ij}$ is a constant depending only on $i,~j~ (1\le i,~j\le k)$. The matrix $A_P=(c_{ij})$ is called the quotient matrix of $G$ with respect $P$.
\end{definition}
\begin{lemma}{\em \cite{CRS}}\label{lem23}
Let $P:V(G)=V_1\cup V_2\cup \cdots \cup V_k$ be an equitable partition of $G$ with quotient matrix  $A_P$. Then $\det(xI-A_P)~|~\det(xI-A(G))$. Furthermore, the largest eigenvalue of $A_P$ is just the spectral radius of $G$.
\end{lemma}

\begin{lemma}{\em\cite{EN}}\label{lem25}
Let $G$ be a graph of size $m$ without isolated vertices. If $G$ is triangle free, then $\rho(G)\le\sqrt{m}$, with equality if and only if $G$ is a complete bipartite graph.
\end{lemma}

\begin{figure}
\centering
\subfigure[]{\includegraphics[width=0.30\linewidth]{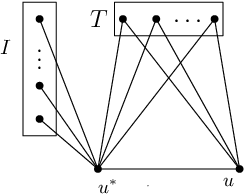}}~~~~~~~~~~
\subfigure[]{\includegraphics[width=0.25\linewidth]{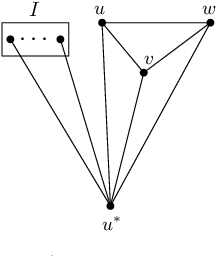}}~~~~
\caption{Graphs (a) $G(m,t)$ and  (b)  $K^m_4$}
\label{fig 6}
\end{figure}

\indent The following lemma shows that the spectral radius of the book graph $S_{\frac{m+3}{2},2}$ is greater than or equal to the spectral radius of $(G(m,t))$.
\begin{lemma}\label{lem24}
Let $G(m,t)$ be the graph as shown in Figure \ref{fig 6}(a) with $m > t+2$. If $t \ge 0$ is even, then $\rho(G(m,t)) \le \rho(S_{\frac{m+3}{2},2})$, and equality holds when $t=0$.
\end{lemma}
\noindent \textbf{Proof.} The quotient matrix of $G(m,t)$ corresponding to the partition $P_3:V(G(m,t))=\{u^\ast\}\cup \{u\}\cup T \cup I$, where $|I|=t$ and $|T|=\frac{m-t-1}{2}$, is
\[
  {A_{P_3}} =
    \bordermatrix{& \{u^\ast\}& \{u\} &  T & I \cr
    \{u^\ast\}& 0 & 1 & \frac{m-t-1}{2} & t \cr
 \{u\} & 1 & 0 & \frac{m-t-1}{2} & 0 \cr
  T & 1 & 1 & 0 & 0 \cr
  I& 1 & 0 & 0 & 0} .\qquad
  \]
Let $f(x) = \det(xI_4 - A_{P_3}) = x^4 - mx^2 - (m - t - 1)x + \frac{t}{2}(m - t - 1)$. The largest root of $S_{\frac{m+3}{2},2}$ satisfies $g(x) = x^3 - mx - (m-1)$. Then $h(x) = f(x) - xg(x) = tx + \frac{t}{2}(m - t - 1)\ge0$ for $x~\textgreater~0$ and $t\ge 0$. This demonstrates that the largest root of $g(x)$ is greater than or equal to the largest root of $f(x)$. Hence the result follows by Lemma \ref{lem23}. \qed

\indent Now, we show that the spectral radius of the book graph is strictly greater than the spectral radius of $K^m_4$.
\begin{lemma}\label{lem28}
Let $K^m_4$ be the graph of size $m$, where $m\ge8$, shown in Figure \ref{fig 6}(b). Then $\rho(K^m_4)~\textless ~ \rho(S_{\frac{m+3}{2},2})$.
\end{lemma}
{\bf Proof.} The vertex set of $K^m_4$ has equitable partition $P_4:V(K^m_4)=\{u^\ast\}\cup R \cup I$, where $R=\{u, v, w\}$ and the quotient matrix with respect to $P_4 $ is
\[
  {A_{P_3}} =
    \bordermatrix{& \{u^\ast\}& R  & I \cr
    \{u^\ast\}& 0 & 3 & m-6 \cr
 R & 1 & 2 & 0 \cr
  I& 1 & 0 & 0 }. \qquad
  \]
  Then $f(x)=\det(xI_3-A_{P_3})=x^3-2x^2-(m-3)x+2(m-6).$ If $\rho_1=\rho(S_{\frac{m+3}{2},2})$, then $\rho_1^2=\rho_1+(m-1)$. Therefore, we have
  \begin{eqnarray*}
  f(\rho_1)&=&\rho_1^2+(m-1)\rho_1-2(\rho_1+m-1)-(m-3)\rho_1+2(m-6)\\
  &=&\rho_1+m-11.
  \end{eqnarray*}
  Given that $m\ge 8$, it follows that $f(\rho_1)>0$. Additionally, the derivative $f^\prime(x)=3x^2-4x-(m-3)>0$ holds for $x\ge \rho_1$. Consequently, by Lemma \ref{lem23}, it can be inferred that $\rho(K^m_4)~\textless ~ \rho(S_{\frac{m+3}{2},2})$.
\qed
In the following lemma, we prove that any vertex of degree one in $G^\ast$ is joined to the extremal vertex $u^\ast$. In particular this shows that $d(v)\ge 2 $ for any $v\in W$.
\begin{lemma} \label{lem26}
Let $G^\ast$ be a graph forbidding the subgraphs $H(3,3)$ and $H(4,3)$, and $u^\ast$ be an extremal vertex in $G^\ast$. Then every pendent vertex in $G^\ast$ is joined to $u^\ast$.
\end{lemma}
\noindent \textbf{Proof.} Suppose there exists a vertex $w \in V(G^\ast) \setminus N(u^\ast)$ with $d(w) = 1$. Let $N(w) =\{v\}$. Consider the graph $G^\prime$ obtained by removing the edge $wv$ from $G^\ast$ and adding the edge $wu^\ast$. The resulting graph $G^\prime$ has $m$ edges and is free of induced subgraphs $H(3,3)$ and $H(4,3)$. Furthermore, $N_{G^\ast}(u^\ast) \subsetneqq N_{G^\prime}(u^\ast)$, and the following inequality holds.
\begin{eqnarray*}
\displaystyle\sum_{\substack{uz \in E(G^\prime)}}x_ux_z& =&\displaystyle\sum_{\substack{uz \in E(G^\ast)}}{x_ux_z } +x_w(x_u^\ast-x_v)   \\
 & \ge& \displaystyle\sum_{\substack{uz \in E(G^\ast)}}{x_ux_z }.
\end{eqnarray*}
According to Lemma \ref{lem21}, this implies that $\rho(G^\prime) > \rho(G^\ast)$. However, this contradicts the definition of $G^\ast$.\qed

\indent In the following lemma we give an upper bound of $e(W)$.
 \begin{lemma}\label{lem27}
Let $G$ be a graph of size $m$ and $W=V(G)\setminus N[u^\ast]$. If $\rho\ge \frac{1+\sqrt{4m-3}}{2}$, then
$$e(W)\le e(N(u^\ast)-|N(u^\ast)\setminus N_0(u^\ast)|+1.$$
\end{lemma}
{\bf Proof.} From Eqs. \ref{eqn(21)} and \ref{eqn(22)}, we can deduce the following expression
\begin{eqnarray*}
(\rho^2-\rho)x_u^\ast &=& d(u^\ast)x_{u^\ast}+\sum_{u\in N(u^\ast)\setminus  N_0(u^\ast)}(d_{N(u^\ast)}(u)-1)x_u+\sum_{w\in W}d_{N(u^\ast)}(w)x_w-\sum_{u\in N_0(u^\ast)}x_u.
\end{eqnarray*}
Given $\rho^2-\rho\ge m-1$, it implies that
\begin{eqnarray*}
(m-1)x_u^\ast &\le& d(u^\ast)x_{u^\ast}+\sum_{u\in N(u^\ast)\setminus  N_0(u^\ast)}(d_{N(u^\ast)}(u)-1)x_u+\sum_{w\in W}d_{N(u^\ast)}(w)x_w-\sum_{u\in N_0(u^\ast)}x_u\\
&\le &\left\{d(u^\ast)+e(N(u^\ast), W)+\sum_{u\in N(u^\ast)\setminus  N_0(u^\ast)}(d_{N(u^\ast)}(u)-1)-\sum_{u\in N_0(u^\ast)}\frac{x_u}{x_u^\ast}\right\}x_u^\ast\\
&=&\left\{m-e(W)+e(N(u^\ast))-|N(u^\ast)\setminus N_0(u^\ast)|-\sum_{u\in N_0(u^\ast)}\frac{x_u}{x_u^\ast}\right\}x_u^\ast.
\end{eqnarray*}
This yields
\begin{eqnarray}\label{eqn23}
e(W) &\le & e(N(u^\ast))-|N(u^\ast)\setminus N_0(u^\ast)|+1-\sum_{u\in N_0(u^\ast)}\frac{x_u}{x_u^\ast}\notag \\
&\le&e(N(u^\ast))-|N(u^\ast)\setminus N_0(u^\ast)|+1.
\end{eqnarray}
This concludes the proof. \qed

\section{Proof of the main theorem}

 Since $G^\ast$ does not contain a subgraph isomorphic to $H(3,3)$, it implies that $G^\ast[N(u^\ast)]$ cannot have two or more independent edges. As a result, $G^\ast[N(u^\ast)]$ can be categorized into one of the following (1) it either consists solely of isolated vertices, or (2) it includes isolated vertices along with a copy of a  star $S_k$, where $k \le \frac{m+1}{2}$, or (3) it comprises isolated vertices and a triangle.\\
\indent The following lemma shows that the only non trivial component of $G^\ast[N(u^\ast)]$ is  a star.

\begin{lemma}\label{lem210}
Let $G^\ast$ be the graph with maximum spectral radius among all $\{H(3,3),H(4,3)\}$-free graphs with odd size $m\ge259$ and let $X^\ast=(x_1,x_2,\cdots,x_n)^T$ be the Perron vector of $G^\ast$, and $x_{u^\ast}=max\{x_i:i\in V(G^\ast)\}$. If $\rho(G^\ast)\ge \frac{1+\sqrt{4m-3}}{2}$, then $G^\ast[N(u^\ast)]$ has exactly one non trivial component $S_k$, where $k\le \frac{m+1}{2}$.
\end {lemma}
\noindent \textbf{Proof.} First, we suppose that $ G^\ast[N(u^\ast)]$ consists of isolated vertices only. Since the star $S_{m+1}$ with $m$ edges does not include $H(3,3)$ and $H(4,3)$ as subgraphs, consequently, due to the extremality of $G^\ast$, it follows that $\rho(G^\ast)\ge \rho(S_{m+1})=\sqrt{m}$. Then, from Eq. \ref{eqn(22)}, we derive the following inequality
\begin{eqnarray*}
(m-d(u^\ast))x_{u^\ast}&=&\sum_{u\in N(u^\ast)\setminus  N_0(u^\ast)}d_{N(u^\ast)}(u)x_u+\sum_{w\in N^2(u^\ast)}d_{N(u^\ast)}(w)x_w\\
&\le&(2e(N(u^\ast))+e((N(u^\ast), N^2(u^\ast)))x_{u^\ast}.
\end{eqnarray*}
This leads to the conclusion that $e(W)\le e(N(u^\ast))$. In other words, $e(W)=0$. According to Lemma \ref{lem26}, there are no pendent vertices in $W$. Consequently, $G^\ast$ is a triangle-free graph. Therefore, by Lemma \ref{lem25}, $G^\ast$ is a complete bipartite graph. Thus, $\rho(G^\ast)=\sqrt{m}~\textless ~\frac{1+\sqrt{4m-3}}{2}$ for $m\ge2$, which is a contradiction.\\
\indent Next, suppose that $G^\ast[N(u^\ast)]$ contains a copy of  $K_3$.  In this case, we observe that $d_{N(u^\ast)}(w)\le 1$. It is clear that $K^m_4$ does not include $H(3,3)$ and $H(4,3)$ as subgraphs. As $S_{m-2}$ is a proper subgraph of $K^m_4$, according to Lemma \ref{lembap}, we must have $\rho(K^m_4)~\textgreater~\rho(S_{m-2})=\sqrt{m-3}$. Moreover, the definition of $G^\ast$ implies that $\rho(G^\ast)~\textgreater~\rho(K^m_4)~\textgreater~\sqrt{m-3}$. Then, based on Eq. \ref{eqn(22)} and using $x_u^\ast\ge x_u$ for any vertex $u$, it can be deduced that
\begin{eqnarray*}
(m-d(u^\ast)-3)x_{u^\ast}&~\textless~&({\rho^\ast}^2-d(u^\ast)x_{u^\ast}\\&\le&(2e(N(u^\ast))+e((N(u^\ast), N^2(u^\ast)))x_{u^\ast}.
\end{eqnarray*}
This implies that $e(W)\le 5$. We claim that $e(W)= 0$. To prove the claim, we consider the following possibilities.\\
{\bf Case 1.} ${\bf 2\le e(W)\le 5}$.\\
\indent Given that $e(W)\le 5$, it follows that $d(w )\le6$ for any $w\in W$ and $d(u)\le8$ for any vertex $u \in N_1(u^\ast)$. Additionally, as $x_u^{\ast}\ge x_u$ for any vertex $u$, we obtain the following inequalities for any vertex $w\in W$ and any vertex $u \in N_1(u^\ast)$.
\begin{eqnarray}\label{eqn(25)}
x_w\le\frac{6}{\rho^\ast}x_u^{\ast}~~\text{and} ~~x_u\le\frac{8}{\rho^\ast}x_u^{\ast}.
\end{eqnarray}
From Eq. \ref{eqn(22)} and inequalities \ref{eqn(25)}, we have\\
\begin{eqnarray*}
{\rho^\ast}^2x_{u^\ast}&=&d(u^\ast)x_{u^\ast}+\displaystyle\sum_{\substack{u\in N(u^\ast)\setminus  N_0(u^\ast)}}d_{N(u^\ast)}(u)x_u+\displaystyle\sum_{\substack{w\in N^2(u^\ast)}}d_{N(u^\ast)}(w)x_w\\
&\le&\Bigg\{d(u^\ast)+\frac{8}{\rho^\ast}\displaystyle\sum_{\substack{u\in N(u^\ast)\setminus  N_0(u^\ast)}}d_{N(u^\ast)}(u)+\frac{6}{\rho^\ast}\displaystyle\sum_{\substack{w\in N^2(u^\ast)}}d_{N(u^\ast)}(w)\Bigg\}x_{u^\ast}\\
&=&\Bigg\{d(u^\ast)+\frac{16}{\rho^\ast}e(N_1(u^\ast)+\frac{6}{\rho^\ast}(e(W,N(u^\ast)))\Bigg\}x_{u^\ast}.
\end{eqnarray*}
Since $\rho^\ast~\textgreater~\sqrt{m-3}\ge16$, from above, we have
\begin{eqnarray}\label{eqn(26)}
{\rho^\ast}^2x_{u^\ast}&\textless&\Bigg\{d(u^\ast)+e(N_1(u^\ast)+\frac{6}{16}(e(W,N(u^\ast)))\Bigg\}x_{u^\ast}\notag\\
&=&\Bigg\{m-e(W)-\frac{10}{16}(e(W,N(u^\ast)))\Bigg\}x_{u^\ast}.
\end{eqnarray}
Recall that $d(w)\ge 2$ for any $w\in W$. Therefore
 \[  e(W,N(u^\ast))\geq \begin{cases}
          2 & e(W)= 2, \\
         1 & e(W)\geq 3 .
       \end{cases}
    \]

Using this in \ref{eqn(26)}, we conclude that $\rho^\ast~\textless~\sqrt{m-3}$, leading to a contradiction.\\

\noindent{\bf Case 2.} ${\bf  e(W)=1}$.\\
\indent In this situation, we observe that $x_w\le\frac{2}{\rho^\ast}x_u^{\ast}$ for any $w \in W$ and $x_u\le\frac{4}{\rho^\ast}x_u^{\ast}$ for any $u \in N_1(u^\ast)$. Following the steps outlined in case 1, we conclude that $\rho^\ast~\textless~\sqrt{m-3}$, contradicting the hypothesis.

Considering both cases 1 and 2, we can deduce that $e(W)=0$. Consequently, the observation that all pendent vertices are joined to $u^\ast$ and $d_{N(u^\ast)}(w)=1$ for any $w \in N^2(u^\ast)$ implies that $W=\emptyset$. Hence, $G^\ast$ is isomorphic to $K^m_4$ when $G^\ast[N(u^\ast)]$ contains a copy of $K_3$. By Lemma \ref{lem28}, we have $\rho(K^m_4)~\textless ~ \rho(S_{\frac{m+3}{2},2})= \frac{1+\sqrt{4m-3}}{2}$, which is again a contradiction. This completes the proof of the lemma.\qed

\noindent{\bf Proof of Theorem 1.1.} As $G^\ast[N(u^\ast)]$ is a tree, by Lemma \ref{lem27}, we have $e(N(u^\ast))-|N(u^\ast)\setminus N_0(u^\ast)| = -1$. Considering that $S_{\frac{m+3}{2},2}$ does not contain $H(3,3)$ and $H(4,3)$ as subgraphs, and by the definition of $G^\ast$, it follows that $\rho(G^\ast)\ge \rho(S_{\frac{m+3}{2},2})=\frac{1+\sqrt{4m-3}}{2}$. According to Lemma \ref{lem27}, this implies that
\begin{eqnarray}\label{eqn24}
e(W) = 0
\end{eqnarray}
and $W=N^2(u^\ast)$. In view of Lemma \ref{lem210}, $ G^\ast[N(u^\ast)]$ includes isolated vertices along with a copy of a  star $S_k$, where $ k \le \frac{m+1}{2}$. Now, we consider the following cases.\\

\noindent {\bf Case 1. $\bf{G^\ast[N(u^\ast)]}$ contains a copy of  ${\bf K_2}$}.\\
\indent  Consider the unique edge $v_1v_2$ in $G^\ast[N_1(u^\ast)]$. In view of Lemma \ref{lem26}, there are no pendent vertices in $W$.  Also $G^\ast$ is restricted from containing the graphs $H(3,3)$ and $H(4,3)$, no vertex in $W$ can have two distinct neighbors in $N_0(u^\ast)$. Therefore, $d_{N(u^\ast)}(w)\leq 3$ for any $w\in W$. \\
\indent Assume that $w\in W$ and $d_{N(u^\ast)}(w)=3$. This implies that $W=\{w\}$ and $w$ is adjacent to $v_1, v_2$ and some vertex $u_1$ in $N_0(u^\ast)$. Let $G_1$ be the graph obtained from $G^\ast$ by deleting the edge $wu_1$ and adding the edge $wu^\ast$. The graph $G_1$ is $\{H(3,3),H(4,3)\}$-free, and $N_{G^\ast}(u^\ast)\subsetneq N_{G_1}(u^\ast)$. Also, $x_{u^\ast}\ge x_{u_1}$. According to Lemma \ref{lem21}, this implies that $\rho(G_1) > \rho(G^\ast)$. However, this contradicts the definition of $G^\ast$. Therefore, $d_{N(u^\ast)}(w)\leq 2$.\\
\indent Suppose that $d_{N(u^\ast)}(w)=2$ for some $w\in W$. If $w$ is adjacent to both $v_1$ and $v_2$, then $N_0^2(u^\ast)\cap N_1^2(u^\ast)=\emptyset$. Consequently, $N^2(u^\ast)=N_1^2(u^\ast)$. Let $w_1, w_2,...,w_s$ be the vertices in $N_1^2(u^\ast)$ that are adjacent to both $v_1$ and $v_2$. Consider $G_2=G^\ast-\{w_iv_1:w_i\in W\}+\{w_iu^\ast:w_i\in W\}$. As previously established by Lemma \ref{lem21}, we derive that $\rho(G_2) > \rho(G^\ast)$, which is not feasible.\\
\indent Now, let each $w \in W$ be adjacent to one vertex in $N_0(u^\ast)$ and one vertex in $N_1(u^\ast)$. Let $w_1, w_2 \in W$ such that $N_{N(u^\ast)}(w_1)=\{v_1,u_1:u_1\in N_0(u^\ast)\}$ and $N_{N(u^\ast)}(w_2)=\{v_2,u_2:u_2\in N_0(u^\ast)\}$.  Without loss of generality, assume that $x_{v_1}\ge x_{v_2}$. Consider $G_3=G^\ast-w_2v_2+w_2v_1$. The graph $G_3$ is $\{H(3,3),H(4,3)\}$-free, and $N_{G^\ast}(v_1) \subsetneq N_{G_3}(v_1)$. Therefore, by Lemma \ref{lem21}, we have $\rho(G_3) > \rho(G^\ast)$, contradicting the definition of $G^\ast$. Thus, every $w\in W$ is adjacent to either $v_1$ or $v_2$ and some vertex $u\in N_0(u^\ast) $. Assume that every $w \in W$ is adjacent to $v_1$. Since $d(w)=2$, so $w$ is adjacent to $v_1$ and some vertex $u_1$ in $N_0(u^\ast)$. Applying the same procedure as before, we obtain the graph $G_4$ from $G^\ast$ by deleting the edge $wu_1$ and adding the edge $wu^\ast$. Clearly, the hypothesis of Lemma \ref{lem21} holds, and it follows that $\rho(G_4) > \rho(G^\ast)$. Thus, we conclude that $W=\emptyset$, and hence $G^\ast$ is isomorphic to $G(m,t)$, where $t=m-3$. By Lemma \ref{lem24}, we have $ \rho(G(m,m-3))~\textless~\rho(S_{\frac{m+3}{2},2})$ and thus contradicting the definition of $G^\ast$. \\
 \begin{figure}[h]
\centering
\subfigure[] {\includegraphics[width=0.35\linewidth]{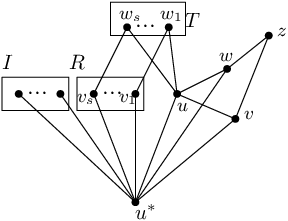}}~~~~~~~~~~~~~
\subfigure[]{\includegraphics[width=0.35\linewidth]{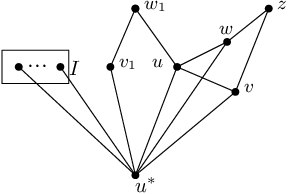}}
\caption{Graphs (a) $H$ and (b) $H_1$}
\label{fig 1}
\end{figure}

\noindent {\bf Case 2. $\bf{G^\ast [N(u^\ast)}]$ contains a copy of star ${\bf S_3}$}.\\
\indent Let $u$, $v$ and $w$ be the vertices of the star component, with $u$ as its center. The following observations apply to any vertex $w\in W$\\
$i$. $w$ cannot be adjacent to any two vertices in $N_0(u^\ast)$.\\
$ii$. $w$ cannot be adjacent to both $u$ and $v$, or both $u$ and $w$.\\
$iii$.  $w$ cannot be adjacent to $v$ (or $w$) and a vertex in $N_0(u^\ast)$.\\
\indent Consider $z \in W$ such that $N_{N(u^\ast)}(z)=\{v, w\}$, and let $T=\{w_1,\cdots,w_s\}$ represent the set of vertices in $W$, and $R=\{v_1,\cdots,v_s\}$ represent the set of vertices in $N(u^\ast)$, satisfying $N_{N(u^\ast)}(w_i)=\{u, v_i\}$ for $i=1,\cdots,s$. Denote this graph by $H$, as depicted in Figure \ref{fig 1}(a). We aim to demonstrate that $\rho(H)~\textless~\rho(H_1)$, where $H_1$ is a graph with $m$ edges and $|T|=1$ as illustrated in Figure \ref{fig 1}(b). To facilitate the analysis, we partition the vertex set of $H$ as $P:\{u^\ast\}\cup \{u\} \cup \{v, w\} \cup R \cup I \cup T \cup \{z\}$, where $|T|=|R|=t$, $|I|=m-3t-7$. The quotient matrix of $H$ concerning the partition $P$ is given by
\[
  {A_P} =
    \bordermatrix{& \{u^\ast\}& \{u\} & \{v, w\} & R & I & T & \{z\} \cr
    \{u^\ast\}   & 0 & 1& 2 & t & m-3t-7&0 &0 \cr
     \{u\} & 1 & 0 & 2 & 0 &0 & t & 0 \cr
      \{v, w\} & 1 & 1 & 0 & 0 & 0 & 0 & 1 \cr
     R & 1 & 0 & 0 & 0 & 0 & 1 & 0 \cr
      I& 1 & 0 & 0 & 0 & 0 & 0 & 0 \cr
       T & 0 & 1 & 0 & 1 & 0 & 0 & 0 \cr
       \{z\}& 0 & 0 & 2 & 0 & 0 & 0 & 0 }.\qquad
     \]
\indent Consider the polynomial $\phi_t(x)=\det(xI_7-A_P)= x^7+(-m+t-1)x^5-4x^4+(mt-2t^2+5m-16t-26)x^3+(-4t+4)x^2+(-2mt+4t^2-4m+30t+26)x$. Consequently, $\phi_1(x) = x^7-mx^5-4x^4+(6m-44)x^3+(-6m+60)x$, and the difference $f(x)=\phi_t(x)-\phi_1(x)=x(t-1)g(x)$, where $g(x)=x^4+(m-2t-18)x^2-4x-2m+4t+34$. To establish $f(x)> 0$, it suffices to show that $g(x)> 0$. Given that $2\le t\le \frac{m-7}{3}$, we have
 \begin{eqnarray*}
 g(x)&\ge& x^4+(m-2(\frac{m-7}{3})-18)x^2-4x-2m+42\\
 &=& \frac{1}{3}(3x^4+(m-40)x^2-12x-6m+126).
 \end {eqnarray*}
 Furthermore, considering that $m\ge 259$, it follows that
\begin{eqnarray*}
g\bigg(\frac{\sqrt{m}}{2}\bigg)=\frac{5}{16}m^2-16m-6\sqrt{m}+126 ~> 0. 
\end {eqnarray*}
Now, for $x~ \textgreater~ \frac{\sqrt{m}}{2}$, the derivative of $g$ is $g^{\prime}(x)=4x^3+\frac{2}{3}(m-40)x-12~\textgreater~ 0 $, and $\frac{\sqrt{m}}{2}\textless \frac{1+\sqrt{4m-3}}{2}$. Consequently, it follows that $f(x)\textgreater0$ for $x~ \textgreater~\frac{1+\sqrt{4m-3}}{2}$, leading to the conclusion that $\rho(H_1)>\rho(H)$. We have the following claim.\\

\begin{figure}[h]
\centering
\subfigure[] {\includegraphics[width=0.30\linewidth]{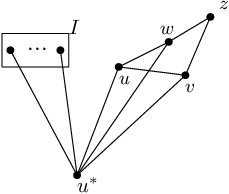}}~~~~~~~~~~~~~~~~~~
\subfigure[]{\includegraphics[width=0.33\linewidth]{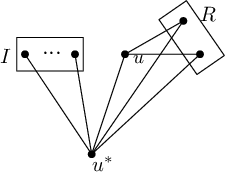}}
\caption{Graphs (a) $H_2$ and (b) $G(m,m-5)$}
\label{fig 3}
\end{figure}

\noindent{\bf Claim.} $\rho(H_1)<\rho(G(m,t)$, where $t=m-5$.  \\
{\textit{ Proof of the claim.}} The quotient matrix of the graph $G(m,m-5)$, as depicted in Figure \ref{fig 3}(b), concerning the partition $P_1:\{u^\ast\}\cup \{u\}\cup R\cup I$ is given by
\[
  {A_{P_1}} =
    \bordermatrix{& \{u^\ast\}& \{u\} &  R & I \cr
    \{u^\ast\}& 0 & 1 &2 & m-5 \cr
 \{u\} & 1 & 0 & 2 & 0 \cr
  R & 1 & 1 & 0 & 0 \cr
  I& 1 & 0 & 0 & 0} .\qquad
  \]
  Consider $\phi_2(x)=\det(xI_4-A_{P_1})= x^4-mx^2-4x+2m-10$. Define
   \begin{eqnarray*}
  \psi_1(x) &=& \frac{1}{x}\phi_1(x)-x^2\phi_2(x) \\& = & x^6-mx^4-4x^3+(6m-44)x^2+(-6m+60)-(x^6-mx^4-4x^3+(2m-10)x^2)\\&=& (4m-34)x^2-6m+60.
  \end{eqnarray*}
  Observe that $\psi^{\prime}_1(x)=2(4m-34)x~\textgreater~0$ and $\psi_1(x)~\textgreater~0$, for $x~\textgreater~\frac{1+\sqrt{4m+3}}{2}$. So it is evident that $ \rho(G(m,m-5))\textgreater \rho(H_1)$. This completes the proof of the claim.\\
  \indent  Consider the case where $|T|=1$, and let $w_1\in T$ be adjacent to both $u$ and $v_1$. In the event that $w$ and $v$ lack a common neighbor in $W$, form the graph $G_1=G^\ast-v_1w_1+u^\ast w_1$. Applying Lemma \ref{lem21} leads to a contradiction.\\
\indent Alternatively, if $T=\emptyset$ and $v$ and $w$ share a common neighbor in $W$, as illustrated by graph $H_2$ in Figure \ref{fig 3}(a).
The quotient matrix of the graph $H_2$ concerning the partition $P_2:\{u^\ast\}\cup \{u\} \cup \{v, w\}\cup \{z\}\cup I$ is given as
\[
 {A_{P_2}} =
    \bordermatrix{& \{u^\ast\}& \{u\} & \{v, w\} &  \{z\} & I \cr
    \{u^\ast\}& 0 & 1 &2 &  0 & m-7 \cr
 \{u\} & 1 & 0 & 2 & 0 & 0 \cr
   \{v, w\} & 1 & 1 & 0 & 1  & 0 \cr
     \{z\} & 0 & 0 & 2 & 0  & 0 \cr
  I& 1 & 0 & 0 & 0 & 0}. \qquad
  \]
\indent  The characteristic polynomial of $A_{P_2}$ is $\phi_3(x)=x^5-mx^3-4x^2+(4m-26)x$. Define $\psi_2(x)=\phi_3(x)-x\phi_2(x)=(2m-16)x~\textgreater~0$ for $x~\textgreater~0$. This implies that $ \rho(G(m,m-5))> \rho(H_2)$.\\
\indent In light of the preceeding discussion, it follows that $G^\ast$ is isomorphic to $G(m,t)$, where $t=m-5$. By Lemma \ref{lem24}, we have $ \rho(G(m,m-5))~\textless~\rho(S_{\frac{m+3}{2},2})$, which contradicts the definition of $G^\ast$.\\

 \noindent{\bf Case 3. $\bf{G^\ast[N(u^\ast)]}$  contains a copy of star ${\bf S_k,~k\ge4}$}.\\
 \indent Consider $u$ as the central vertex of the star $S_k$. With $\rho(G^\ast)\ge \frac{1+\sqrt{4m-3}}{2}$, it follows that $e(W)=0$. Furthermore, as $G^\ast$ is free of induced subgraphs $H(3,3)$ and $H(4,3)$, it ensures that no two vertices in $N_0(u^\ast)$ or $N_1(u^\ast)$ share a common neighbor in $W$. If $w \in N^2_0(u^\ast) \cap N^2_1(u^\ast)$, then $|N_{N(u^\ast)}(w)|=2$ and $w$ must be adjacent to $u$. Let $N_{N(u^\ast)}(w)=\{u, v: ~v\in N_0(u^\ast)\}$. Define $G_1=G^\ast-vw+u^\ast w$. Applying Lemma \ref{lem21} in this context leads to a contradiction. Consequently, $W=\emptyset$ and by Lemma \ref{lem24}, it follows that $G^\ast$ is isomorphic to $S_{\frac{m+3}{2},2}$. This completes the proof of the theorem.  \qed

\section{Conclusion}

\indent In this paper, we have established that the graph denoted by $S_{\frac{m+3}{2},2}$ has the maximum spectral radius within the class of \{$H(3,3),H(4,3)$\}-free graphs with odd size greater than or equal to 259. For cases where $m$ is less than 259, we identify Lemma 3.1 as a crucial technical obstacle. To investigate and to overcome this obstacle, extending our proof to the case where $m<259$ would be an interesting problem for further research. Moreover, it's important to highlight that Theorem 1.1 establishes the proof for the case of odd $m$. This naturally raises the question: What is the maximum spectral radius among $\{H(3,3),H(4,3)\}$-free graphs when the size is even? This question invites further investigation and we leave it as an open problem.\\

\noindent{\bf Conflict of interest.} The authors declare that they have no conflict of interest.\\

\noindent{\bf Data Availibility} Data sharing is not applicable to this article as no datasets were generated or analyzed
during the current study.\\

\end{document}